\documentclass{article}
\author{Igor~Bayak}
\title{On the relation of the monomial group \\
with other algebraic structures}
\begin{document}
\maketitle
\noindent {\bf Summary:} {\it It is shown
in what way monomial group connects Abelian
group $Z^{n}$ and total linear group $GL(n)$.
It is shown that any subgroup of Abelian group
$Z^{n}$ induces subgroup of monomial group
$S_{2}\wr S_{n}$, which in its turn induces
corresponding subgroup of $GL_{n}(R)$.}
\section{Monomial group as a generator of total
linear group}

In combinatorial constructions of this section
we use language of set and map theory therefore,
succeeding its notions and notations, let
$I=\{1,\ldots ,n\}$;  $s:I\to I$~--- is bijection,
$s's:I\to I:s's(i)=s'\left(s(i)\right)$~--- is the
composition of bijetions. Then set of all $n!$
bijections $S=\{s\}$ with respect to product
$S\times S\to S:(s',s)\to s's$ put together
permutation group (symmetric group) of the order $n$.
At the same time, set of all even bijections put
together alternating group $S^{+}$, which in general
case is generated by 3--cycles, namely,
$S_{n}^{+}:=\left\langle
\{(i,j,k)\}_{I}\right\rangle$ where $n\ge 3$ and
$i\ne j\ne k$.

Let also $A=\{\pm 1,\ldots ,\pm n\}$ and
$p:I\to A$~--- is such an injection that $\vert p
\vert :I\to I$~--- is bijection, moreover
$p'p:I\to A:p'p(i)=\textrm{sgn}\, p(i)\cdot p'(
\vert p(i)\vert)$~--- is the composition of these
injections. Then set of all $2^{n}n!$ modulo
bijective injections $P=\{p\}$ with respect to
product $P\times P\to P:(p',p)\to p'p$ put
together monomial group described by the wreath
product $P_{n}:=S_{2}\wr S_{n}$ and intuitively
imagined as group of arrowy permutations
(arrowy group) of the order $n$. Next let arrowy
permutation is called as arrowy transposition
if it is a simple elementary permutation or
inversion of one arrow, in other words, element
of monomial group is transposition if it is
an elementary permutation in $S_{n}$--component
or a rearrangement in one from $S_{2}$--component
of the wreath product. Then, like the simple
permutations, parity of arrowy permutation $p$
is defined by parity of number of arrowy
transpositions for passage from identity
permutation to $p$, which is invariant relative
to the choise of composition of transpositions for
passage to $p$. In this connection, the set of all
even arrowy permutations put together group
$P^{+}$, which in general case is generated
by rearrangements with inversion, i.e. by generators
$(j,-k):j\to(-k),k\to j$ where $j\ne k$, so we
have generation $P_{n}^{+}:=\left\langle\{(j,-k)\}_{I}
\right\rangle$.

However, along with the notion of parity of
permutations (arrowy permutations) it should be
considered the notion of parity of arrangements
(arrowy arrangements), for which it is sufficient
to accept that elementary (transpositional)
arrangement onto the place $I^{\ast}$ (where
$I^{\ast}\subset I$) corresponds to either
elementary rearrangement into the set $I^{\ast}$
or one--element replacement from complementary
to $I^{\ast}$ set  $I\backslash I^{\ast}$
(the notion of elementary arrowy arrangement must
be extended by addition of one--element inversion
onto the place $I^{\ast}$). Then parity of arrangements
$s^{\ast}:=s(I^{\ast})$ and $p^{\ast}:=p(I^{\ast})$
is defined by factored at $Z_{2}$ numbers
$\sigma(s^{\ast})$, $\sigma(p^{\ast})$ of elementary
arrangements, that are required for passage to
$s^{\ast}$ and  $p^{\ast}$ respectively.

The notion of arrangement parity allows to form
new combinatorial groups. In fact, let
$I_{1}=\{1,\ldots ,r\}$, $I_{2}=\{r+1,\ldots ,n\}$,
$s^{1}:=s(I_{1})$, $s^{2}:=s(I_{2})$, then we can
present $s$ as a component permutation, consisting of
arrangements $s^{1}$, $s^{2}$ and next define component
parity given by number $\sigma(s^{1})+\sigma(s^{2})$.
Hence, set of all the even--component permutations are gathered
in group $S_{r,n-r}^{+}$, which in general case is
generated by 3--cycles, acting into subsets $I_{1},I_{2}$
and by arbitrary 2--cycles, acting between them,
namely, $$S_{r,n-r}^{+}:=
\left\langle\{(i,j,k)\}_{I_{1}},(l,m),
\{(i,j,k)\}_{I_{2}}\right\rangle$$ where
$l\in I_{1}$ and $m\in I_{2}$.

In its turn, set of all the even--component arrowy
permutations put together group $P_{r,n-r}^{+}$,
which in general case is generated by generators
$(i,-j)$, acting into subsets $I_{1},I_{2}$, and
by arbitrary pair generator $\pm (l,m)$, acting
between them, where $+(l,m)$~--- is a simple
rearrangement and $-(l,m)$~--- is an inverse
rearrangement, i.e. $-(l,m):l\to (-m),m\to (-l)$.
Hence, we have generation:
$$P_{r,n-r}^{+}:=\left\langle\{(i,-j)\}_{I_{1}},\pm
(l,m),\{(i,-j)\}_{I_{2}}\right\rangle$$ where
$l\in I_{1}$ and $m\in I_{2}$.

Further, before to enter upon the matrix realization
of monomial group we define the determinant of row
subset of square matrix $A=\left (a_{i,j}\right
)_{i,j\in I}$ (i.e. of rectangular matrix
$A(I^{\ast}):=\left(a_{i,j}\right)
_{i\in I^{\ast},j\in I}$)
as a sum of products of elements indexed
by arrangements of columns on rows and multipled at the
sign of parity of this arrangement. In other words,
we assume that $$\det{A(I^{\ast})}:=
\sum_{\{s^{\ast}\}}\prod_{i\in I^{\ast}}
(-1)^{\sigma(s^{\ast})}\cdot a_{is^{\ast}(i)}$$

Let now $P(I):=\left\{\left(p_{ij}\right)
_{I}\right\}$~--- is an exhaustive set of such
square matrices in which every column and every
row have only one nonzero element, moreover if
$p_{ij}\ne 0$, then $p_{ij}\in\{\pm 1\}$. Then
we have linear presentation $P(I)$ of arrowy
group $P_{n}$ given by formula
$p_{ij}=\textrm{sgn}\, p(i)\cdot 1$ 
for every $j=\vert p(i)\vert$ and $p_{ij}=0$ for
every $j\ne\vert p(i)\vert$. Besides, matrix groups
$$\left\{\left (p_{ij}\right )_{I}\vert\det{P(I)}
=+1\right\}$$ $$\left\{\left (p_{ij}\right )_{I}
\vert\det{P(I^{\ast})}\cdot\det{P(I\backslash
I^{\ast})}=+1\right\}$$  are linear presentations
of even and even--component arrowy groups respectively.
At the same time, if given set $J=\{1,\ldots ,m\}$
where $m\le n$ and family $I_{J}:=\left\{I_{j}\right\}_{J}$
where $\cup_{m}I_{j}=I$, $\cap_{m}I_{j}=\oslash$,
$\textrm{card}\, I_{j}=n_{j}$, $\sum_{m}n_{j}=n$, moreover
subset $I_{j+1}$ is filled by sequential sampling from
$I$ after filling $I_{j}$ and $I_{1}=\{1,\ldots ,n_{1}\}$,
then matrix group $$\left\{\left(p_{ij}\right
)_{I}\vert\det{P(I_{1})}\cdots\det{P(I_{m})}=+1
\right\}$$
gives linear presentation of arbitrary subgroup
$P^{+}_{n_{1},\ldots ,n_{m}}$ of even--component
arrowy permutations.

At last, we establish the connection of
monomial group with total linear group. Let
$RP_{n}$~--- is algebra over $R$, which is
generated as the linear span of group $P_{n}$
represented as a subset of matrix algebra $ML_{n}(R)$,
i.e. $RP_{n}=R\langle P(I)\rangle$; $GP_{n}$~---
is multiplicative group of algebra $RP_{n}$;
$R^{\ast}$~--- is multiplicative group of field $R$;
$\textrm{Aut}\, RP_{n}:=\textrm{Int}\, GP_{n}
\simeq GP_{n}/R^{\ast}$~--- is group
of automorphisms of this algebra. Then we get
underlying statement:$$GP_{n}=GL_{n}(R)$$

In fact, since group $P_{n}$ is realized
by set of transitional matrices $P(I)$ from which
we always can choose $n^{2}$ linearly independent
matrices, then it includes some basis of algebra
$ML_{n}(R)$, hence the linear span of this group coincides
with $ML_{n}(R)$, i.e. $RP_{n}=ML_{n}(R)$, and therefore
$GP_{n}=GL_{n}(R)$. In consequence of underlying statement
we have equality $\textrm{Aut}\, RP_{n}=PGL_{n}(R)$ 
where $PGL_{n}(R)$~--- is the total projective group.

Let now $GL_{ij}:=\textrm{diag}\left
[1,\ldots ,GL(2)^{i,j}_{i,j},\ldots ,1\right]_{n}$
where $i\ne j$, moreover upper indexes point to
rows numbers and lower indexes point to columns
numbers onto intersection of which is disposed
group $GL(2)$. Then as additional source of
procedure of generation of linear groups will
be statement:
$$GL(n)=\left\langle\left\{GL_{i,j}\right\}
_{I}\right\rangle=\left\langle\left\{GL
_{i,i+1}\right\}_{I}\right\rangle$$

Really, since $P_{n}=\left\langle\left
\{P^{2}_{i,j}\right\}_{I}\right\rangle=\left
\langle\left\{P^{2}_{i,i+1}\right\}_{I}\right\rangle$
where $P^{2}_{i,j}$~--- is arrowy subgroup
of the order $n$ presented by the arrowy permutations
onto place $i,j$ and isomorphic to group  $P_{2}$,
then $GP_{n}=\left\langle\left\{GP^{2}
_{i,j}\right\}_{I}\right\rangle$, but
$GP^{2}_{i,j}=GL_{i,j}$, hence the statement is proved.

Further let there be given linear presentations
of two arrowy groups
$P^{+}_{2}=\left\langle\left
(\begin{array}{cc}0 & 1 \\ -1 & 0
\end{array}\right)\right\rangle$ and
$P^{+}_{1,1}=\left\langle\pm\left
(\begin{array}{cc}0 & 1 \\ 1 & 0
\end{array}\right)\right\rangle$, then after
elementary calculations we obtain that
$$\textrm{Aut}\, RP^{+}_{2}=\left\{\left
(\begin{array}{cc}\cos{x}&\sin{x}\\
-\sin{x}&\cos{x}\end{array}\right)
,x\in R\right\}=SO(2)$$ and
$$\textrm{Aut}\, RP^{+}_{1,1}=
\left\{\pm\left(\begin{array}{cc}
\cosh{x}&\sinh{x} \\ \sinh{x}&\cosh{x}
\end{array}\right),x\in R\right\}=SO(1,1)$$
and denoting
$$SO(2)_{i,j}:=
\textrm{diag}\left[1,\ldots ,SO(2)^{i,j}_{i,j}
,\ldots ,1\right]_{n}$$
and
$$SO(1,1)_{i,j}:=\textrm{diag}\left[1,\ldots 
,SO(1,1)^{i,j}_{i,j},\ldots ,1\right]_{n}$$
obtain generation of special orthogonal group
$$\textrm{Aut}\, RP^{+}_{n}=\left\langle
\left\{SO(2)_{i,i+1}\right\}_{I}\right
\rangle=SO(n)$$ and special pseudo--orthogonal
group $$\textrm{Aut}\, RP^{+}_{r,n-r}
=\left\langle\left\{SO(2)_{i,i+1}\right\}
_{I_{1}},SO(1,1)_{r,r+1},\left\{SO(2)
_{i,i+1}\right\}_{I_{2}}\right\rangle=SO(r,n-r)$$
However, generalization of special orthogonal
group is possible, including special case of
pseudo--orthogonal group, namely,
$$\textrm{Aut}\, RP^{+}_{n_{1},\ldots,n_{m}}
=SO(I_{J}):=SO(n_{1},\ldots,n_{m})=$$
$$\left\langle\left\{\left\{SO(2)
_{i,i+1}\right\}_{I_{j}}\right\}_{J},
\left\{SO(1,1)_{n_{j},n_{j+1}}\right\}
_{J}\right\rangle$$

In the whole, it is clear that all
possible subgroups of monomial group
generate all possible linear groups
(among them~--- unitary, symplectic, exeptional
and so on), therefore the problem of classification
and construction of linear groups coincides
with the one for subgroups of monomial group.

\section{Monomial group as group of automorphisms of 
$Z$--modulus $Z^{n}$}

Let $Z$~--- is ring of integers; $Z^{n}$~--- is
$n$--dimensional $Z$--modulus and
$e_{I}:=\left(e_{i}\right)_{I}$~--- is its basis,
then $\textrm{Aut}\, Z^{n}:=\textrm{Int}\,
GL_{n}(Z)\simeq P_{n}$, since from the equality
$GL_{n}(Z)/\{\pm 1\}=P(I)$ it follows that
$P(I)$~--- is group of inner automorphisms of
total integer linear group $GL_{n}(Z)$,
i.e. of group of nonsingular integral--valued
matrices. Further we will look for
group of such automorphisms of modulus $Z^{n}$,
which correspond to topology of the quotient
$\left\{-e_{I},o,e_{I}\right\}/\sim$, i.e. to
partition of the set consisting of base elements,
inverse base elements and the point of reference
at equivalence classes. But, since the partition
of this base set must be induced by factor group
$Z^{n}/H^{n}$, then it should be said on
automorphisms of modulus $Z^{n}$, corresponding to
topology of factor group. For example, for the factor
group $Z/Z$ we have one class
$[-e\sim e\sim o]$ and for the factor group $Z/mZ$,
where $m\ge 2$, we have partition on two classes
$[-e\sim e],[o]$, so topological type of the
generalized factor group $Z/H$ may be
presented by groups $Z/Z$ and $Z/2Z$.
It is characteristic also that any factor
group $Z^{n}/H^{n}$ induces at least
equivalence relation $-e_{i}\sim e_{i}$
for every $i\in I$.

Turning to automorphisms of base set,
we note that map $-I\times I\to -P(I)\times P(I)$
realizes coordinated with automorphisms $Z^{n}$
rearrangement of elememts of set
$\left\{-e_{I},e_{I}\right\}$ by formula
$$-e_{I}\times e_{I}\to -e'_{I}\times e'_{I}:e'_{i}=e_{p(i)}:=
\textrm{sgn}\, p(i)\cdot e_{\vert p(i)\vert}$$
Moreover, reversible maps:
$$-e_{i}\times e_{i}\to e_{i}\times -e_{i}$$
$$-e_{j},-e_{k}\times e_{j},e_{k}\to -e_{k}
,-e_{j}\times e_{k},e_{j}$$
$$-e_{j},-e_{k}\times e_{j},e_{k}\to e_{k}
,e_{j}\times -e_{k},-e_{j}$$
realize elementary reversible rearrangements of elements of set
$\left\{-e_{I},e_{I}\right\}$, coordinated with
arrowy inversion, simple rearrangement and
inverse rearrangement respectively. If we
present now quotient $\left\{-e_{I},o,e_{I}
\right\}/\sim$ as $2n$--link graph glued in
equivalence points, then it will be clear that
rebuilding of this factored graph may be not
solved, namely, it is impossible to realize
rearrangement of two links, the ends of which are
glued with each other but not glued with the point of
reference, where 2--link element of graph it
should be concidered as 2--dimensional, i.e.
planar graph.

Thus, starting from topological arguments,
we postulate the prohibition of reversible rearrangement
of arbitrary pair elements of set
$\left\{-e_{I},e_{I}\right\}$ if they are
equivalent in some class, different from the
class of the point of reference. Hence, we have
restriction of the group of automorphisms $Z^{n}$
to some subgroup of monomial group.
So if pairs $\left(-e_{i},e_{i}\right)$ of
base elements are grouped in $m$ classes,
none from which is not equivalent to the point of
reference, then group of topological automorphisms
$Z^{n}$ corresponds to subgroup
$P^{+}_{n_{1},\ldots ,n_{m}}$. Really, only
generators of even--component arrowy group
(i.e. rearrangement with inversion) exclude
permutations, corresponding to all the elementary
reversible rearrangements in each set
$\left\{-e_{I_{j}},e_{I_{j}}\right\}$, but include
simple and inverse rearrangements between each
$I_{j}$ and  $I_{k}$, that corresponds to
choosen factorization of set
$\left\{-e_{I},o,e_{I}\right\}$.

We present also convenient method of
construction of arbitrary factor group.
Let $\vert z\vert\bmod{2}:Z\to Z/2Z:
z'\equiv z\pmod{2}:
\vert z\vert\bmod{2}=r$, where $r$~--- is
remainder from division $z\colon 2$,
taken without the sign of $z$, then we have
field $Z_{2}=\langle\{0,1\};\oplus ,
\cdot,0,1\rangle$, where $x\oplus y=\vert
x+y\vert\bmod{2}$. At the same time, if
$Z^{n}_{2}:=\prod_{n}Z_{2}$,
then $Z^{n}_{2}$~--- is Abelian group with
component--wize addition and moreover
$Z^{n}_{2}$~--- is linear space over $Z_{2}$
and $e_{I}:=(e_{i})_{I}$, where
$e_{i}=(0_{(1)},\ldots ,1_{(i)},\ldots ,0_{(n)})$,~---
is its proper base. Let also
$H^{n}_{2}$~--- is the notation of the generalized
subgroup of group $Z^{n}_{2}$ and $H(I)$~---
is the concrete subgroup, consisting of all
possible even sum of elements of proper base,
namely,
$H(I)=\left\langle\left\{e_{j}\oplus e_{k}\right\}_{I}
\right\rangle$, which divides $Z^{n}_{2}$ at
two coset, so $Z^{n}_{2}/H(I)\approx Z_{2}$.
Let us assume $H(I_{J}):=\coprod_{m}H(I_{j})$,
then $Z^{n}_{2}/H(I_{J})\approx Z^{m}_{2}$.
At the same time, since we have homomorphism
$z:H^{n}_{2}\to H^{n}:(0\to 2Z,1\to 2Z+1)$, then
any factor group $Z^{n}_{2}/H^{n}_{2}$ is
isomorphic to some factor group $Z^{n}/H^{n}$,
which in its turn induces the group of topological
automofphisms of $Z^{n}$, in particular
$\textrm{Aut}\, Z^{n}/z\ast H(I_{J})=
P^{+}_{n_{1},\ldots ,n_{m}}$ and
$\textrm{Aut}\, Z^{n}/z\ast Z^{n}_{2}=
\textrm{Aut}\, Z^{n}/Z^{n}=P_{n}$.

Further, keeping in mind that any quotient
$\left\{R^{n}/H^{n}\right\}$ may be identified
with some compact space, topology of which is
defined by factor group $Z^{n}/H^{n}$, we obtain
identification with $n$--sphere
$\left\{R^{n}/z\ast H(I)\right\}\approx S^{n}$
or with more of generality with toroidal space
$\left\{R^{n}/z\ast H(I_{J})\right\}\approx 
S^{n_{1}}\times\cdots\times S^{n_{m}}$.
At the same time, since $RP_{n}=\textrm{End}\, R^{n}$,
then $\textrm{End}\, S^{n_{1}}\times\cdots
\times S^{n_{m}}\approx RP^{+}_
{n_{1},\ldots ,n_{m}}$ and
$\textrm{Aut}\, S^{n_{1}}\times\cdots\times
S^{n_{m}}\approx SO(n_{1},\ldots ,n_{m})$, hence the
rigidity of Euclidean spaces stipulated by
certain type of functional of scalar product
is the result of narrowing of algebra of endomorphisms
and group of automorphisms of the space $R^{n}$,
arising from its factorization. In particular,
geometry of proper Euclidean space is associated
with factorization of $n$--dimensional
arithmetical space into $n$--dimensional
sphere. On the other hand, projective
space $RP^{n}\approx\left\{R^{n}/Z^{n}\right\}$
does not possess topological rigidity
and therefore $\textrm{End}\, RP^{n}=RP_{n}$
and $\textrm{Aut}\, RP^{n}=PGL_{n}(R)$.

Thus, monomial group by own subgroups provides
conformity between subgroups of Abelian group $Z^{n}$
and subgroups of total linear group $GL_{n}(R)$,
moreover any Abelian subgroup is associated with
some topologically compact space, group of automorphisms
of which coincides with corresponding linear subgroup.

\end{document}